\documentclass{amsart}
\usepackage{latexsym}
\usepackage{amsmath, amssymb, amsfonts, amsxtra}
\usepackage{amsthm}

\author{Richard Oberlin}
\email{oberlin@math.wisc.edu}
\address{University of Wisconsin-Madison, Mathematics Department, 480 Lincoln Dr, Madison WI 53706}

\title{Bounds for Kakeya-type maximal operators associated with $k$-planes}
\date{}

\newcommand {\leb}{\mathcal{L}}
\newcommand {\gra}{\mathcal{G}}
\newcommand {\ort}{\mathcal{O}}
\newcommand {\rea}{\mathbb{R}}
\newcommand {\sph}{\mathbb{S}}
\newcommand {\ma}{\mathcal{M}}

\newcommand {\proj}{\mathrm{proj}}
\newcommand {\spa}{\mathop{\mathrm{span}}}
\newcommand {\na}{\mathcal{N}}

\newcommand {\Holder} {H\"{o}lder}
\newcommand {\kcrit} {k_{cr}}

\newcommand {\composed} {\circ}

\newtheorem{lem}{Lemma}[section]
\newtheorem{thm}{Theorem}[section]
\newtheorem{cor}{Corollary}[section]
\newtheorem{conj}{Conjecture}[section]
\newtheorem{prop}{Proposition}[section]

\begin{document}

\begin{abstract}
A $(d,k)$ set is a subset of $\rea^d$ containing a translate of
every $k$-dimensional plane. Bourgain showed that for $k \geq \kcrit(d)$,
where $\kcrit(d)$ solves 
$2^{\kcrit-1}+\kcrit = d$, every $(d,k)$ set has positive Lebesgue measure. 
We give a short proof of this result which allows for an improved $L^p$
estimate of the corresponding maximal operator, and which demonstrates 
that a lower value of $\kcrit$ could be obtained if improved mixed-norm estimates for the $x$-ray transform were known.
\end{abstract}

\maketitle

\section{Introduction}

A measurable set $E \subset \rea^d$ is said to be a $(d,k)$ set if it contains 
a translate of every $k$-dimensional plane in $\rea^d$. Once the definition
is given, the question of the minimum size of a $(d,k)$ set arises. This question 
has been extensively studied for the case $k=1$, the Kakeya sets. It is known 
that there exist Kakeya sets of measure zero, and these are called Besicovitch sets. It is 
conjectured that all Besicovitch sets have Hausdorff dimension $d$. For $k \geq 2$, it is 
conjectured that $(d,k)$ sets must have positive measure, i.e. that there are no $(d,k)$ 
Besicovitch sets. These size estimates are related to $L^p$ bounds on two
maximal operators which we define below. 

Let $G(d,k)$ denote the Grassmannian manifold 
of $k$-dimensional linear subspaces of $\rea^d$. For $L \in G(d,k)$ we define
\[
\na^k [f](L) = \sup_{x \in \rea^d} \int_{x + L} f(y) dy
\]
where we will only consider functions $f$ supported on the unit ball $B(0,1) 
\subset \rea^d$.

A limiting and rescaling argument shows that if $\na^k$ is bounded for some $p < \infty$ 
from $L^p(\rea^d)$ to 
$L^1(G(d,k))$, then $(d,k)$ sets must have positive measure. By testing $\na^k$ on the 
characteristic function of $B(0,\delta)$, $\chi_{B(0,\delta)}$, 
one sees that such a bound may only hold for $p \geq \frac{d}{k}$. 
For $ L $ in $ G(d,k) $ and $ a \in \rea^d $ define the $ \delta $ plate
centered at $a$, $ L_\delta(a) $, to be the $\delta$ neighborhood in $ \rea^d $
of the intersection of $ B(a,\frac{1}{2}) $ 
with $ L + a $. Fixing $L$, considering $\na^k \chi_{L_\delta(0)}$, and using the fact that 
the dimension of $G(d,k)$ is $k(d-k),$ we see that a bound into $L^q(G(d,k))$ can only hold for $q \leq kp$.
This leads to the following conjecture, where the case $k=1$ is excluded due to the existence of 
Besicovitch sets.
\begin{conj} \label{conjbdna}
For
$2 \leq k < d, p > \frac{d}{k}, 1 \leq q \leq k p$
\[
\|\na^k f \|_{L^q(G(d,k))} \lesssim  \|f\|_{L^p(\rea^d)}.
\] 
\end{conj}
It is also useful to consider a generalization of the Kakeya 
maximal operator, defined for $L \in G(d,k)$ by
\[
\ma^k_{\delta}[f](L) = \sup_{a \in \rea^d} \frac{1}{\leb^d(L_\delta(a))} 
\int_{L_\delta(a)} f(y) dy
\]
where $\leb^d$ denotes Lebesgue measure on $\rea^d$.
Using an argument analogous to that in Lemma 2.15 of \cite{bg}, one may see that 
a bound 
\begin{equation} \label{bound}
\|\ma^k_\delta f \|_{L^1(G(d,k))} \lesssim \delta^{\frac{-\alpha}{p}} \|f\|_{L^p(\rea^d)}
\end{equation}
where $ \alpha > 0 $ and $ p < \infty $, implies that the Hausdorff dimension of any $ (d,k) $ 
set is at least $ d - \alpha $.
Considering $\ma^k_\delta \chi_{B(0,\delta)}$ and $\ma^k_\delta \chi_{L_\delta(0)}$, we formulate
\begin{conj} \label{conjbdma}
For $k \geq 1, p < \frac{d}{k}, q \leq (d-k)p'$
\[
\|\ma^k_\delta f \|_{L^q(G(d,k))} \lesssim  \delta^{k-\frac{d}{p}} \|f\|_{L^p(\rea^d)}.
\] 
\end{conj}
In \cite{fa1} Falconer showed that, for any $\epsilon > 0$, $\na^k$ is bounded from $L^{\frac{d}{k}+ \epsilon}(\rea^d)$ 
to $L^1(G(d,k))$ when 
$k > \frac{d}{2}$. Later, in \cite{bg}, Bourgain used a Kakeya maximal operator bound
combined with an $L^2$ estimate of the $x$-ray transform to show that $\na^k$ is bounded 
from $L^p(\rea^d)$ to $L^p(G(d,k))$ for $(d,k,p) = (4,2,2 + \epsilon)$ and $(d,k,p) = (7,3,3 + \epsilon)$. 
He then showed, using a recursive metric entropy estimate, that   
for $ d \leq 2^{k-1} + k $, $\na^k$ is bounded for a large unspecified $p$. 
Substituting in the proof Katz and Tao's more recent bound for the Kakeya maximal operator from \cite{kt}
\begin{equation} \label{ktkakbound}
\| \ma^1_\delta f \|_{L^{n+\frac{3}{4}}(G(n,1))} \lesssim \delta^{-\left(\frac{3(n-1)}{4n+3} + \epsilon\right)}
 \|f\|_{L^{\frac{4n+3}{7}}(\rea^n)}
\end{equation}
one now sees that this holds for $k > \kcrit(d)$
where 
\begin{equation} \label{kcritdef}
\kcrit(d) \text{\ solves\ } d = \frac{7}{3} 2^{\kcrit-2} + \kcrit.
\end{equation}

By \Holder 's inequality, the following is true for any $k$-plate $L_\delta$ and positive $f$
\[
\int_{L_\delta}f \ dx \lesssim
\delta^{\frac{d-k}{r'}} 
\left( \int_{L^\perp} \left(\int_{L+y} f(x)\ d\leb^k(x) \right)^r d\leb^{d-k}(y)   \right)^{\frac{1}{r}}.
\]
Combining this with the $L^p \rightarrow L^q(L^r)$ bounds for the $k$-plane transform which were proven by Christ in 
Theorem A of \cite{ch}, we see that Conjecture \ref{conjbdma} holds with 
$ p \leq \frac{d+1}{k+1} $. Except for a factor of $\delta^{-\epsilon}$, the same bound
for $ \ma^k_\delta $ was proven with $k=2$ by Alvarez in \cite{al} using a geometric-combinatorial
``bush''-type argument. Alvarez also used a ``hairbrush'' argument to show that $(d,2)$ sets have Minkowski dimension at least $\frac{2d+3}{3}$. More recently, Mitsis proved a similar maximal operator bound in \cite{mi2} and showed that $(d,2)$ sets have Hausdorff dimension at least $\frac{2d+3}{3}$ in \cite{mi1}. In \cite{bu}, Bueti extended these dimension estimates, in the context of finite fields, to $(d,k)$ sets, showing that $(d,k)$ sets in $F^d$ have dimension at least $\frac{k(d+1)+1}{k+1}$. In \cite{rogers}, Rogers gave estimates for the Hausdorff
dimension of sets which   
contain planes in directions corresponding to certain curved submanifolds of $G(4,2)$.

Our main result is the following.
\begin{thm} \label{boundonnak}
Suppose $4 \leq k < d$ and $k > \kcrit(d)$, where $\kcrit(d)$ is defined in (\ref{kcritdef}). Then 
\begin{equation} \label{finalnakbound}
\|\na^k f\|_{L^{p}(G(d,k))} \lesssim 
\|f\|_{L^{p}(\rea^d)}
\end{equation}
for $f$ supported on the unit ball and $p \geq \frac{d-1}{2}.$
If, additionally, we have $k-j > \kcrit(d-j)$ for some integer $j$ in $[1,k-4]$, 
then we may take $p \geq \frac{d-1}{2+j}.$ 
\end{thm} 

For $k < \kcrit(d)$, we do not have a bound for $\na^k$, however our technique yields certain bounds for $\ma^k_\delta$. 

\begin{thm} \label{mabounds}
\[
\|\ma^k_\delta f\|_{L^q(G(d,k))} \lesssim \delta^{-\frac{\alpha}{p}}
\|f\|_{L^p(\rea^d)}
\]
when
\begin{equation} \label{sharpp}
k \geq 2, \ \alpha = d - kp + \epsilon, \ p = \frac{d}{k+\frac{3}{4}}, \ q \leq (d-k)\left(\frac{4(d-(k-1))}{7}\right)'
\end{equation}
or
\begin{equation} \label{nonL2ma}
k \geq 2, \ \alpha = \frac{3(d-k)}{7(2^{k-1})}+ \epsilon, \ p=\frac{d+1}{2}, \ q=d+1
\end{equation}
or
\begin{equation} \label{L2ma}
3 \leq k \leq \kcrit(d), \ \alpha = \frac{3(d-k)}{7(2^{k-2})} - 1 + \epsilon, \ p = q = \frac{d}{2}
\end{equation}
where $\epsilon > 0$ may be taken arbitrarily small.
\end{thm}

In  (\ref{sharpp}) we have an optimal value for $p$ relative to $\alpha$, 
but a non-optimal value for $q$. In (\ref{nonL2ma}) and (\ref{L2ma}) we 
have improved values of $\alpha$ at the cost of a non-optimal $p$. 
For the ``non-borderline'' $k$, specifically when $k+1 < \kcrit(d+1)$, (\ref{nonL2ma}) gives a smaller value of $\alpha$ than (\ref{L2ma}).  

The number $p = \frac{d-1}{2+j}$ in Theorem \ref{boundonnak} and the number
$p = \frac{d}{k+\frac{3}{4}}$ in Theorem \ref{mabounds} are approximate and may be slightly improved through careful numerology. Also, in (\ref{L2ma}) we may take $k = 2$, but a slightly higher value of $p$ and $q$ is then required. 

We prove (\ref{sharpp}) and (\ref{nonL2ma}) in Section \ref{armob}
through a recursive maximal operator bound which is derived using 
Drury and Christ's bounds for the $x$-ray transform and which is 
inspired by Bourgain's recursive metric entropy estimates. 
This recursive maximal operator bound is a slight improvement of the result in 
\cite{ro}, which will remain unpublished, and the new bound comes with a vastly simplified proof afforded by the explicit use of 
the $x$-ray transform.
Additionally our argument
reveals that with certain adjustments of $p$ and $q$, the number $2$ in the definition of $\kcrit(d)$ and 
in the definition of $\alpha$ in (\ref{nonL2ma}) and (\ref{L2ma}) 
may be replaced by the ratio $\frac{\tilde{r}}{\tilde{p}}$ if the $x$-ray transform 
is known to be bounded, for certain values of $n$, from $L^{p_n}(\rea^n)$ to $L^{q_n}_{\sph^{n-1}}(L^{r_n}_{\rea^{n-1}})$ for any $r_n,p_n,q_n$ satisfying $\frac{r_n}{p_n} = \frac{\tilde{r}}{\tilde{p}}$. 

We prove  
(\ref{L2ma}) and Theorem \ref{boundonnak} in Section \ref{fourthsection}. There, we combine 
(\ref{sharpp}) and (\ref{nonL2ma}) with the $L^2$ method which Bourgain 
used to give bounds for $\na^k$ when $(d,k) = (4,2)$ or $(7,3)$. 

From (\ref{nonL2ma}) and (\ref{L2ma}) we see that, for $k \geq 2$, the 
Hausdorff dimension of any $(d,k)$ set is at least 
\[
\mathrm{min}\left(d, \max\left(d- \frac{3(d-k)}{7(2^{k-2})} + 1, d- \frac{3(d-k)}{7(2^{k-1})}\right)\right).
\] 
When $(d-k) < 7$, it is preferable to start with Wolff's $L^{\frac{n+2}{2}}$ bound for the Kakeya maximal operator from \cite{wo2}, instead of (\ref{ktkakbound}). A similar procedure then gives the lower bound
\[
\mathrm{min}\left(d, \max\left(d - \frac{d-k-1}{2^{k-1}} + 1, d - \frac{d-k-1}{2^k}\right)\right)
\]
for the Hausdorff dimension of a $(d,k)$ set.

It should be noted that the dimension estimates provided by applying 
(\ref{nonL2ma}) and it's Wolff-variant 
are also a direct consequence of the metric entropy estimates in \cite{bg}. 
However, to the best of the author's knowledge 
they have not previously appeared in the literature, even without the improvement obtained from \cite{wo} and \cite{kt}.

\subsection*{Acknowledgements}
I would like to thank my advisor Andreas Seeger for mathematical 
guidance and for his suggestion of the topics considered in this article.
I would also like to thank Dan Oberlin for carefully reading several 
of the previous drafts.

\section{A recursive maximal operator bound} \label{armob}

We start with the definition of the measure we will use on $G(d,k)$.
Fix any $L \in G(d,k)$. For a Borel subset $F$ of $G(d,k)$ let
\[
\gra^{(d,k)}(F)=\ort(\{\theta \in O(d): \theta(L) \in F\})
\]
where $\ort$ is normalized Haar measure of the orthogonal group on $\rea^d$, $O(d)$.  
By the transitivity of the action of $O(d)$ 
on $G(d,k)$ and the invariance of $\ort$, it is clear that the definition is 
independent of the choice of $L$. Also note that $\gra^{(d,k)}$ is invariant under the 
action of $O(d)$. By the uniqueness of uniformly-distributed measures (see \cite{mt}, pages 44-53),
$\gra^{(d,k)}$ is the unique normalized Radon measure on $G(d,k)$ invariant under $O(d)$.

It will be necessary to use an alternate formulation of $\gra^{(d,k)}$. For each $ \xi $ in 
$ \sph^{d-1} $ let $ T_\xi:\xi^\perp \rightarrow \rea^{d-1} $ be an orthogonal linear 
transformation. Then $ T_\xi^{-1} $ identifies $ G(d-1,k-1) $ with the $ k-1 $ dimensional 
subspaces of $ \xi^\perp $. Now, define $T:\sph^{d-1}\times G(d-1,k-1) \rightarrow G(d,k)$ by
\[
T(\xi,M)=\spa(\xi,T_\xi^{-1}(M)).
\]
Choosing $T_\xi$ continuously on the upper and lower hemispheres of $ \sph^{d-1} $, 
$T^{-1}$ identifies the Borel subsets of $G(d,k)$ with the completion of the Borel
subsets of $\sph^{d-1} \times G(d-1,k-1)$. Under this identification, by uniqueness of 
rotation invariant measure, we have
\begin{equation} \label{graproduct}
\gra^{(d,k)}(F) = \sigma^{d-1} \times \gra^{(d-1,k-1)}(T^{-1}(F))
\end{equation}
where $\sigma^{d-1}$ denotes normalized surface measure on the unit sphere.

For a function $f$ on $\rea^d$, $\xi \in \sph^{d-1}$, and $y \in \xi^{\perp},$ the $x$-ray transform of $f$ is defined  
\[
f_\xi(y) = \int_{\rea} f(y + t\xi)\ dt.
\]
It is conjectured that the $x$-ray transform is bounded from $L^p(\rea^d)$
to $L^q_{\sph^{d-1}}(L^r_{\rea^{d-1}})$ when $p,q,r$ satisfy
\begin{eqnarray} 
\nonumber r &<& \infty
\\ \label{rpcondition} p &=& \frac{rd}{d + r - 1}
\\ \nonumber
q &\leq& r'd.
\end{eqnarray}
This was shown to hold in \cite{dr} for $p < \frac{d+1}{2}$ and in \cite{ch} for
$p=\frac{d+1}{2}$. Also, see \cite{wo} and \cite{lt} for certain improvements.

In the following proposition we exploit the fact that 
$r > p$ when $r \neq 1$ in (\ref{rpcondition}), i.e. that the $x$-ray 
transform is $L^p$-improving. 

\begin{prop} \label{maxopone}
Suppose that $p \leq d+1$ and $k \geq 2$. Then a bound 
\[
\|\ma_\delta^{k-1}f \|_{L^q(G(d-1,k-1))} \lesssim \delta^{-\frac{\alpha}{p}} 
\|f\|_{L^p(\rea^{d-1})}
\]
for all $f \in L^{p}(\rea^{d-1})$
implies the bound 
\[
\|\ma_\delta^{k}f \|_{L^{\tilde{q}}(G(d,k))} \lesssim \delta^{-\frac{\tilde{\alpha}}{\tilde{p}}} 
\|f\|_{L^{\tilde{p}}(\rea^{d})}
\]
for all $f \in L^{\tilde{p}}(\rea^{d})$ with
\[
\tilde{p} = p\frac{d}{d+p-1}, \ \ \ \tilde{\alpha} = \alpha\ \frac{\tilde{p}}{p} = \alpha\ \frac{d}{d + p - 1}, \ \ \ \text{and}\ \ \ \tilde{q} = \min(q,dp').
\]
\end{prop}

\begin{proof} 
Without loss of generality, we assume that $f$ is positive. 
Let $L \in G(d,k)$ and suppose that $L = \spa(\xi, T_\xi^{-1}(M))$ where $M \in G(d-1,k-1)$. 
Let $a_L \in \rea^d$ and let $a_M = T_\xi(\proj_{\xi^{\perp}}(a_L))$, where 
$\proj$ denotes orthogonal projection.  
Then 
\begin{eqnarray*}
\int_{L_\delta(a_L)} f(y)\ dy \leq \int_{M_\delta(a_M)} \int_{\rea} f(T_\xi^{-1}(x) + t \xi)\ dt\ dx
\\ = \int_{M_\delta(a_M)} f_\xi(T_\xi^{-1}(x))\ dx
\end{eqnarray*}
where $L_{\delta}(a_L)$ and $M_{\delta}(a_M)$ are $k$ and $k-1$ plates respectively. Noting that $d-k$ = $(d-1)-(k-1)$, it follows that  
\[
\ma_\delta^{k}[f](L) \lesssim \ma_\delta^{k-1}[f_\xi \composed T_\xi^{-1}](M).
\]

By (\ref{graproduct}), H\"{o}lder's inequality, and our hypothesized bound, we now have
\begin{eqnarray*}
\|\ma_\delta^{k}[f]\|_{L^{\tilde{q}}(G(d,k))} 
\lesssim \left( \int_{\sph^{d-1}} \int_{G(d-1,k-1)} 
\ma_\delta^{k-1}[f_\xi \composed T_\xi^{-1}](M)^{\tilde{q}}\ dM\ d\xi
 \right)^{\frac{1}{\tilde{q}}}
\\\lesssim \left( \int_{\sph^{d-1}} \left(\int_{G(d-1,k-1)} 
\ma_\delta^{k-1}[f_\xi \composed T_\xi^{-1}](M)^{q}\ dM\right)^{\frac{\tilde{q}}{q}}\ d\xi
 \right)^{\frac{1}{\tilde{q}}} 
\\ \lesssim \delta^{-\frac{\alpha}{p}} 
\left( \int_{\sph^{d-1}} \left( \int_{\rea^{d-1}} (f_\xi \composed T_\xi^{-1}(x))^{p}\ dx\right)^{\frac{\tilde{q}}{p}} d\xi \right)^{\frac{1}{\tilde{q}}}
\\ =  \delta^{-\frac{\alpha}{p}} 
\left( \int_{\sph^{d-1}} \left( \int_{\xi^{\perp}} f_\xi(x)^{p}\ dx\right)^{\frac{\tilde{q}}{p}} d\xi \right)^{\frac{1}{\tilde{q}}}.
\end{eqnarray*}

Finally, by our restrictions on $p$ and $\tilde{q}$, we may apply Drury and Christ's 
bound for the $x$-ray transform, obtaining 
\begin{eqnarray*}
\left( \int_{\sph^{d-1}} \left( \int_{\xi^{\perp}} f_\xi(x)^{p}\ dx\right)^{\frac{\tilde{q}}{p}} d\xi \right)^{\frac{1}{\tilde{q}}} \lesssim
\|f\|_{L^{\tilde{p}}(\rea^{d})}
\end{eqnarray*}
when $\tilde{p} = \frac{pd}{d+p-1}.$

\end{proof}

One should note that if $\alpha = (d-1) - (k-1)p$, then 
$\tilde{\alpha} = d - k \tilde{p}$. Hence, except for a non-optimal
$\tilde{q}$, Proposition \ref{maxopone} yields the conjectured bound on $L^{\tilde{p}}(\rea^d)$ 
when applied to the conjectured bound on $L^p(\rea^{d-1})$. 

\begin{proof}[Proof of (\ref{sharpp})]
Observing that if \begin{equation} \label{pestimate}
p = \frac{d-1}{m} \ \ \ \text{then} \ \ \ \tilde{p} = \frac{(d+1)-1}{m+1},
\end{equation}
we start from the bound 
\begin{equation} \label{weakerktkakbound}
\| \ma^1_\delta f \|_{L^{(n-1)\left(\frac{4n}{7}\right)'}} \lesssim \delta^{-(\frac{3}{4} + \epsilon)}
\|f\|_{L^{\frac{4n}{7}}(\rea^n)}
\end{equation}
with $n = d-(k-1)$, which is weaker but more convenient for numerology than 
(\ref{ktkakbound}). Since (\ref{weakerktkakbound}) satisfies the left side of (\ref{pestimate}) with $m = \frac{7}{4}$ and $d=n+1$, we obtain (\ref{sharpp}) after $k-1$ iterations  of Proposition \ref{maxopone}. 
\end{proof}
For a larger improvement in $\alpha$, one may interpolate the 
known $L^p$ bound for $\ma^{k-1}_\delta$ with the trivial $L^\infty$ bound 
and apply Proposition \ref{maxopone} to the resulting $L^{d+1}$ bound. 
This allows us to use the maximum value, $2$, of $\frac{r}{p}$ permitted
by Drury and Christ's bound, and yields the following corollary.

\begin{cor} \label{halfalpha}
Under the assumptions of Proposition \ref{maxopone}, we may also take
$\tilde{p}=\frac{d+1}{2}$, $\tilde{\alpha}=\frac{\alpha}{2}$, and $\tilde{q} = \min(\frac{(d+1)q}{p},(d+1))$.
\end{cor}

Due to the interpolation, Corollary \ref{halfalpha} cannot yield a bound for
which $\alpha$ is sharp with respect to $p$ as in Conjecture \ref{conjbdma}.

\begin{proof}[Proof of (\ref{nonL2ma})]
Starting from (\ref{ktkakbound}) with $n=d-(k-1)$, we iteratively apply
Corollary \ref{halfalpha} $(k-1)$ times to obtain (\ref{nonL2ma}).
\end{proof} 

We would like to point out that
the proof of Proposition \ref{maxopone} and Corollary \ref{halfalpha} is similar in spirit to Bourgain's recursive metric
entropy estimate in the sense that a more efficient version of his technique, namely the proof of Proposition 3.1 in \cite{ro}, 
could be used to derive the localized non-endpoint version of the $L^{\frac{d+1}{2}} \rightarrow L^{d+1}$ $x$-ray transform 
bound. The idea of expressing an average over a $k$-plane as the average over a $k-1$-plane of the $x$-ray transform and then ``unraveling'' the integration 
over $G(d,k)$ into a product integral over $\sph^{d-1}$ and $G(d-1,k-1)$ is also due to Bourgain, as he used it in Propositions 3.3 and 3.20 of \cite{bg}. There, he gave bounds for 
$\na^k$ with  $(d,k)=(4,2)$ and $(d,k)=(7,3)$. We state a generalization of this result below, omitting a few details from the proof, as it is essentially the 
same as in \cite{bg}. 
 
\section{The $L^2$ method} \label{fourthsection}
Reducing $\alpha$ by a factor of two, as in Corollary \ref{halfalpha}, is not a substantial gain
for small $\alpha$. By using an $L^2$ estimate of the $x$-ray transform which takes advantage of 
cancellation, instead of the $L^{\frac{d+1}{2}}$ bound, we may take $\tilde{\alpha} = \alpha - 1$ when $\alpha \geq 1$ and obtain a bound for $\na^k$ when 
$\alpha < 1$. 

\begin{prop} \label{bourgainrecursive}
Suppose $k,p \geq 2$ and that a bound for $\ma^{k-1}_\delta$ on $L^p(\rea^{d-1})$ of the form 
\begin{equation} \label{assumedboundbr} 
\|\ma^{k-1}_\delta f\|_{L^p(G(d-1,k-1))} \lesssim \delta^{-\frac{\alpha}{p}} 
\|f\|_{{L^p}(\rea^{d-1})}
\end{equation}
is known. Then if $\alpha \geq 1$ we have the bound
\begin{equation} \label{alphageq1}
\|\ma^k_\delta f\|_{L^{p}(G(d,k))} \lesssim 
\delta^{-\frac{\alpha - 1}{p}}\|f\|_{L^{p}(\rea^d)}
\end{equation}
for $f \in L^p(\rea^d)$.
If $\alpha < 1$ we have the bound
\begin{equation} \label{negativealpha}
\|\na^k f\|_{L^{p}(G(d,k))} \lesssim 
\|f\|_{L^{p}(\rea^d)}
\end{equation}
for $f \in L^{p}(\rea^d)$ supported on $B(0,1)$.
\end{prop}

\begin{proof}[Proof of Theorem \ref{boundonnak}] 
We start from the bound (\ref{nonL2ma})
with $d_0 = d - 2 -j$ and $k_0 = k - 2 - j$. 
This gives
\begin{equation} \label{d0conditions}
\alpha_0 = \frac{3(d-k)}{7\cdot 2^{k-3-j}} + \epsilon, \ p_0 = \frac{d_0 + 1}{2}, \ \text{\ and\ } q_0 = d_0 + 1. 
\end{equation}
The condition $k-j > \kcrit(d-j)$
ensures that $\alpha_0 < 2$, and so no further improvement in $\alpha$ 
is necessary. Thus, 
we use our $j$ ``spare'' iterations to improve $p$. 
We note that, in Proposition \ref{maxopone}, when $m \leq d$,
\begin{equation} \label{L2pcond}
p \leq \frac{d}{m} \ \ \text{implies that\ }\ \tilde{p} \leq \frac{d+1}{m+1}.
\end{equation}
Since $p_0$ satisfies the left inequality in (\ref{L2pcond}) with $m = 2$ and $d = d_0 + 1$, we see that we may take
\[
p_1 = \frac{d_1 + 1}{3},\ q_1 = d_0 + 1,\ \text{\ and\ } \alpha_1 = \alpha_0,
\]
where $d_1 = d_0 + 1 = d-2 - (j-1)$ and $k_1 = k_0 - 1 = k - 2 - (j-1)$.
Above, we ignore the improvement in $\alpha$ and, through interpolation, we ignore some slight additional improvement in $p$. After $j-1$ further iterations, 
we have
\begin{equation} \label{twomoreiterations}
p_j = \frac{d_j + 1}{2+j},\ q_j = d_0 + 1,\ \text{\ and\ } \alpha_j = \alpha_0,
\end{equation}  
where $d_j = d-2$ and $k_j = k-2$.
Applying (\ref{alphageq1}) to (\ref{twomoreiterations}), and then applying (\ref{negativealpha}) to the result, we obtain (\ref{finalnakbound}).
\end{proof}
  
\begin{proof}[Proof of (\ref{L2ma})]
We obtain (\ref{L2ma}) by starting from (\ref{nonL2ma}) with $d_0 = d-1$, and $k_0 = k-1$ (In the case $k=2$, we would simply start from (\ref{ktkakbound})). We then apply (\ref{alphageq1}) once. 
\end{proof}

To derive Proposition \ref{bourgainrecursive}, we use the following estimate. Below, $\hat{f}$ denotes the Fourier transform of $f$.

\begin{lem} \label{cancellation}
Suppose $\hat{f} \equiv 0$ on $B(0,R)$. Then
\[
\|f_\xi(y)\|_{L^2_{\xi,y}(\sph^{d-1} \times \rea^{d-1})} \lesssim R^{-\frac{1}{2}}\|f\|_{L^2(\rea^d)}.
\]
\end{lem}

The above lemma was proven in \cite{bg}, but we give a different
proof which yields a slightly stronger result. 

\begin{lem} \label{radontransform}
For $d \geq 3$
\[
\|f_\xi(y)\|_{L^2_{\xi,y}(\sph^{d-1} \times \rea^{d-1})} = C_d \|f\|_{\dot{H}^{-\frac{1}{2}}(\rea^d)}
\]
where $C_d$ is a fixed constant depending only on $d$ and $\dot{H}$ denotes the homogeneous $L^2$ Sobolev space.
\end{lem}

\begin{proof}
Applying Plancherel's theorem to the partial Fourier transform in the $ \xi^\perp $ direction, we have for every $\xi \in \sph^{d-1}$
\[
\int_{\xi^\perp} |f_\xi(x)|^2 dx =
\int_{\xi^\perp} \left| \hat{f}(\zeta) \right|^2 d\zeta 
\]
where $\hat{f}$ denotes the full Fourier transform of $f$.
Then
\[
\int_{\sph^{d-1}} \int_{\xi^\perp} |f_\xi(x)|^2 dx\ d\xi  
= \int_{\sph^{d-1}} \int_{\rea^{d-1}} |\hat{f}\composed T_\xi^{-1}(\zeta)|^2 d\zeta\ d\xi  
\]
where $T_{\xi}^{-1}$ is defined as in Section \ref{armob}.
Using polar coordinates in the $\zeta$ variable gives
\[
\int_{\sph^{d-1}} \int_{\rea^{d-1}} |\hat{f}\composed T_\xi^{-1}(\zeta)|^2 d\zeta\ d\xi = C
\int_{\sph^{d-1}} \int_{\sph^{d-2}} \int_{\rea} |\hat{f}\composed T_\xi^{-1}(\omega r)|^2 r^{d-2}\ dr\ d\omega\ d\xi.
\]
By the uniqueness of rotation invariant
measures on $\sph^{d-1}$, we have for every $g$ 
\[
\int_{\sph^{d-1}} \int_{\sph^{d-2}} g(T_{\xi}^{-1}(\omega)) d\omega \ d\xi
=
\widetilde{C} \int_{\sph^{d-1}} g(\xi) d\xi. 
\]
Then, since 
$T_{\xi}^{-1}(\omega r) = r T_{\xi}^{-1}(\omega)$
\begin{eqnarray*}
\int_{\sph^{d-1}} \int_{\sph^{d-2}} \int_{\rea} |\hat{f}\composed T_\xi^{-1}(\omega r)|^2 r^{d-2}\ dr\ d\omega\ d\xi
&=& \widetilde{C} \int_{\sph^{d-1}} \int_{\rea}|\hat{f}(r\xi)|^2 r^{d-2}\ dr\  d\xi \\
&=& \overline{C} \int_{\rea^{d}} |\hat{f}(\zeta)|^2 |\zeta|^{-1} \ d\zeta.
\end{eqnarray*}
\end{proof}

Let $f$ be a nonnegative function supported on the unit ball in $\rea^{d}$. To apply Lemma \ref{cancellation}, we use a Littlewood-Paley decomposition, writing 
\[
f = \sum_{j = 0}^{\infty} f_j
\]
where $f_j = f * \phi_j$, $\hat{\phi}_0 = \chi_{B(0,1)}$, and $\hat{\phi}_{j} = \chi_{B(0,2^j)} - \chi_{B(0,2^{j-1})}$ for $j > 0$.

Since $f$ is supported on the unit ball, we may switch the order of integration between convolution and the $x$-ray transform to obtain
\[
\|(f_j)_{\xi}(y)\|_{L^\infty_{\xi,y}} \lesssim \|f\|_{L^{\infty}(\rea^d)}
\]
uniformly in $j$. Hence, interpolation with Lemma \ref{cancellation} gives
\begin{equation} \label{L2interpolant}
\|(f_j)_{\xi}(y)\|_{L^p_{\xi,y}}  \lesssim (2^{-j})^{\frac{1}{p}} \|f\|_{L^p(\rea^d)}
\end{equation}
for any $p \geq 2$. 
Following the proof of Proposition \ref{maxopone}, we observe that for 
$L = \spa(\xi,T_{\xi}^{-1}(M))$ we have 
\begin{equation} \label{macomposed}
\ma^k_{\delta}[f](L) \lesssim \ma^{k-1}_{\delta}[f_{\xi}\composed T_{\xi}^{-1}](M).
\end{equation}
Approximating $\chi_{M_{\delta}}$ by a version with compact Fourier-support 
and estimating the Schwartz-tails, one sees that 
\[
\ma^{k-1}_{\delta}(g) \lesssim \ma^{k-1}_{\delta}(|\tilde{g}|)
\]
for nonnegative functions $g$, and any function $\tilde{g}$ which satisfies $\hat{\tilde{g}}=\hat{g}$ on $B(0,\frac{1}{\delta})$.
Thus, we obtain 
\begin{equation} \label{L2step1} 
\ma^{k-1}_{\delta}[f_{\xi}\composed T_{\xi}^{-1}](M) \lesssim \sum_{j=0}^{|\log(\delta)|+1}
\ma^{k-1}_{\delta}[|(f_j)_{\xi}\composed T_{\xi}^{-1}|](M).
\end{equation}
Another Schwartz-tail estimate shows that for each $j$
\begin{equation} \label{L2step2}
\ma^{k-1}_{\delta}[|(f_j)_{\xi} \composed T_{\xi}^{-1}|](M) \lesssim 
\ma^{k-1}_{2^{-j}}[|(f_j)_{\xi} \composed T_{\xi}^{-1}|](M).
\end{equation}

Integrating over $G(d,k)$ and combining the bounds (\ref{assumedboundbr}) and (\ref{L2interpolant}) as in the proof of Proposition \ref{maxopone}, we obtain
\[
\|\ma^{k}_{\delta}f\|_{L^{p}(G(d,k))} \lesssim \sum_{j=0}^{|\log{\delta}|+1}(2^{j})^{\frac{\alpha-1}{p}}\|f\|_{L^p(\rea^d)} \lesssim \delta^{-\frac{\alpha-1}{p}} \|f\|_{L^{p}(\rea^d)}
\]
from (\ref{macomposed}), (\ref{L2step1}), and (\ref{L2step2}), when $\alpha \geq 1$.

Similarly, we have
\[
\na^k[f](L) \lesssim \na^{k-1}[f_{\xi}\composed T_{\xi}^{-1}](M)
\]
and
\[
\na^{k-1}[|(f_j)_{\xi} \composed T_{\xi}^{-1}|](M) \lesssim 
\ma^{k-1}_{2^{-j}}[|(f_j)_{\xi} \composed T_{\xi}^{-1}|](M),
\]
giving
\[
\|\na^{k}f\|_{L^{p}(G(d,k))} \lesssim \sum_{j=0}^{\infty} (2^{j})^{\frac{\alpha-1}{p}}\|f\|_{L^p(\rea^d)} \lesssim \|f\|_{L^{p}(\rea^d)}
\]
when $\alpha < 1$.

\end{document}